\documentclass[11pt]{amsart}

\usepackage{amsmath, amssymb}
\usepackage{amscd}
\usepackage{verbatim}
\newtheorem{definition}{Definition}[section]
\newtheorem{theorem}[definition]{Theorem}

\theoremstyle{definition}
\newtheorem{remark}[definition]{Remark}

\newtheorem{unitary representation}[definition]{Unitary representation}
\newtheorem{building of discrete dual}[definition]{Building of discrete dual}
\newtheorem{left regular representation}[definition]{Left regular representation}
\newtheorem{action of a discrete quantum group}[definition]{Action of a discrete quantum group}
\newtheorem{weak expectation property}[definition]{Weak expectation property}
\newtheorem{invariant mean}[definition]{Invariant mean}

\newcommand{\clA}{\mathcal{A}}
\newcommand{\hG}{\widehat{\mathbb{G}}}

\newcommand{\cG}{\mathbb{G}}
\newcommand{\hL}{\widehat{\lambda}}
\newcommand{\hW}{\widehat{\omega}}
\newcommand{\clH}{\mathcal{H}}

\newcommand{\hD}{\widehat{\Delta}}

\newcommand{\AG}{\clA\rtimes_{\alpha } \widehat{\mathbb{G}}}
\newcommand{\bbW}{\mathbb{W}}

\newcommand{\id}{\mbox{id}}
\newcommand{\irrG}{\mbox{Irr}\mathbb{G}}

\begin{document}


\title[WEP of DQG Algebras and Crossed Products]{Weak Expectations of Discrete Quantum Group Algebras and Crossed Products}

\author[A.~Bhattacharjee]{Arnab Bhattacharjee}
\address{Department of Mathematics, IISER Bhopal, MP 462066, India \newline (\textit{Present Address}) Mathematical Institute of Charles University, Prague, Czech Republic 18675}
\email{bhattacharjeea@karlin.mff.cuni.cz}

\author[A.~Bhattacharya]{Angshuman Bhattacharya}
\address{Department of Mathematics, IISER Bhopal, MP 462066, India}
\email{angshu@iiserb.ac.in}



\begin{abstract}
In this article we study analogues of the weak expectation property of discrete group C*-algebras and their crossed products, in the \textit{discrete quantum group} setting, i.e., \textit{discrete quantum group} C*-algebras and crossed products of C*-algebras with \textit{amenable discrete quantum groups}.
\end{abstract}

\maketitle


In this article the notation $\otimes$ denotes the minimum tensor product of C*-algebras.

\section{Introduction}

The \textit{weak expectation property} was introduced by Lance \cite{L} in his study of nuclearity of (discrete) group C*-algebras. A \textit{weak expectation} of a represented C*-algebra $A\subset B(H)$ is a unital completely positive map $$\Phi: B(H) \rightarrow A^{\prime\prime}$$ such that $\Phi|_A = \mbox{id}_A$, where $A^{\prime\prime}$ denotes the double commutant of $A$ in $B(H)$. A C*-algebra is said to have the \textit{weak expectation property} if $A$ admits a weak expectation for every \textit{faithful} representation $A\subset B(K)$.

We recall a few well-known facts about the weak expectation property. For a discrete group $\Gamma$, the weak expectation property of the reduced group C*-algebra $C_r^*(\Gamma)$ is equivalent to the \textit{amenability} of the group $\Gamma$, see \cite{L, BO}. Also, in general, the weak expectation property has poor permanence properties, unlike those of nuclearity, exactness and the \textit{quotient weak expectation property} (see \cite{Kir93, Oz04, BO} for definitions, examples and permanence results). However, as a particular instance, where permanence does hold: it was shown in \cite{BF} that, under an amenable action of a discrete group $\Gamma$ on a C*-algebra $A$, the C*-crossed product $A\rtimes \Gamma$ has the weak expectation property \textit{if and only if} so has $A$.

Let $\cG$ be a compact quantum group and $\widehat{\cG}$ denote the discrete quantum dual of $\cG$. The purpose of this article is to investigate the following in the context of the \textit{discrete quantum group} $\widehat{\cG}$:
\begin{enumerate}
\item[1.] Equivalence of existence of a weak expectation $\Phi: B(L^2(\cG))\rightarrow C^*_r(\hG)^{\prime\prime}$ and the amenability (or co-amenability) of $\widehat{\cG}$ (or $\cG$).  
\item[2.] Permanence property for weak expectations of the crossed product $A\rtimes_\alpha \widehat{\cG}$ and $A$, for amenable $\widehat{\cG}$. 
\end{enumerate}

While an answer to (1.) in full generality is yet to be discovered (owing to the non-traciality of the Haar functional of a general compact quantum group), a close analogue to the discrete group case is obtained by the existence of a special weak expectation, which we introduce in Definition \ref{qwe} and call \textit{quantum weak expectation}. Our definition of quantum weak expectation is motivated from \cite{SV14}. See also \cite{CN}. We have the following theorem in this regard:
\begin{theorem} \label{dqgwe}
Let $\cG$ be a compact quantum group and let $\hG$ denote the discrete quantum dual of $\cG$. Then $C^*_r(\hG)$ admits a quantum weak expectation if and only if $\hG$ is amenable.
\end{theorem}

The permanence question in (2.) is completely settled and is analogous to the discrete group case. To achieve this, we use Skalski and Zacharias' finite rank approximation results of the crossed product algebra by a discrete quantum group \cite{SZ10} and give an  explicit constructive proof. In this case we have the following result:  
\begin{theorem} \label{cr} 
Let $\alpha$ be a faithful action of an amenable discrete quantum group $\hG$ on a unital $C^*$-algebra $\clA$, then $\clA\rtimes_{\alpha }\hG$ has the weak expectation property if and only if $\clA$ does.
\end{theorem}

\begin{remark}
From a natural perspective, Theorem \ref{cr} may be regarded as a complement to the permanence results in \cite[Theorem 4.4]{SZ10} in the discrete \textit{quantum} group setting, just as \cite[Theorem(s) 2.1, 3.1]{BF} complements \cite[Theorem 4.3.4]{BO} in the discrete group setting.
\end{remark}

This article is organized as follows: in section 2 we briefly recall the necessary definitions and constructions required for our purpose. In section 3 we present the proof(s) of Theorem(s) \ref{dqgwe} and \ref{cr}. 

\section{Preliminaries}

We briefly recall some definitions and notations used in this article. This section, by no means, aim to be a detailed recollection of facts. The readers are directed to the references cited herein for a comprehensive exposition to the subject matter.
  
\subsection*{Compact quantum groups}\cite{Wor87b,Wor95} A C*-algebraic \textit{compact quantum group} is a pair $\cG= (A,\Delta)$, where $A$ is a unital $C^*$-algebra, and $\Delta :A\rightarrow A\otimes A$ is a unital $*$-homomorphism such that:
\begin{itemize}
\item $\Delta$ is coassociative: $(\Delta \otimes \id)\circ \Delta= (\id\otimes \Delta)\circ \Delta$
\item $\overline{\mbox{span}}((1\otimes A) \Delta (A))= A\otimes A= \overline{\mbox{span}}((A\otimes 1) \Delta (A))$
\end{itemize}

We denote the equivalence classes of irreducible representations of the compact quantum group $\cG$ by $\irrG$.

\subsection*{Discrete quantum groups}\cite{PWor,VV06,Wor95} The \textit{Pontryagin dual} of a compact quantum group $\cG$ is a \textit{discrete quantum group} denoted by $\hG$. The discrete dual $\hG$ (C*-algebraic version) is a pair $\hG=(c_0(\hG),\hD)$ or (vonNeumann algebraic version) $\hG=(\ell^\infty(\hG),\hD)$ where:
$$c_0(\hG)=\bigoplus_{\alpha\in \irrG} B(H_\alpha), \qquad \ell^\infty(\hG)=\prod_{\alpha\in \irrG} B(H_\alpha)$$
and $\hD$ is the comultiplication on $c_0(\hG)$ or $\ell^\infty(\hG)$. If $\hG=\Gamma$ (discrete group) then we have the respective algebras as $c_0(\Gamma)$ and $\ell^\infty(\Gamma)$. 

\subsection*{Regular representation} The \textit{multiplicative unitary} $\mathbb{W}\in B(L^2(\cG)\otimes L^2(\cG))$ associated with the discrete quantum group $\hG$ is such that: $$\hD (x)=\mathbb{W}^*(1\otimes x)\mathbb{W},$$ where $\hD$ is the \textit{co-multiplication} of $\hG$. Let $\hW \in \ell^1(\hG):=\ell^\infty(\hG)_*$. The \textit{left regular representation} $\hL:\ell^1(\hG)\rightarrow B(L^2(\cG))$ of $\mathbb{\hG}$ is defined by: $$\hL(\hW)=(\hW \otimes \mbox{id})\mathbb{W}.$$The \textit{reduced} discrete quantum group algebra is defined as: $$C^*_r(\hG):=\overline{\hL(\ell^1(\hG)).}^{||\cdot||_{B(L^2(\cG))}}$$ By \textit{Peter-Weyl theory}, we have $C^*_r(\hG)=A$ where $\cG= (A,\Delta)$. Further, for compact quantum groups, we have $\ell^\infty(\hG)\subset B(L^2(\cG))$ as a von-Neumann subalgebra (see \cite{Wor87b, Wor95, FSS17, B} for details).

\subsection*{Amenability}\cite{BMT2, BT03, Tomatsu06, B} An \textit{invariant mean} $m$ on $\hG$ is a state $m\in \ell^{\infty}(\hG)^*$ satisfying $m(\hW \otimes \mbox{id}) \hD = \hW (1)m$. A discrete quantum group $\hG$ is said to be \textit{amenable} if there exists an invariant mean on $\ell^{\infty}(\hG)$. This definition of amenability of $\hG$ is one of the many equivalent ones in the literature (see the references mentioned above for equivalent definitions and other details).

\subsection*{Actions and crossed products}\cite{VV06, B, SZ10} A (left) action of a discrete quantum group $\hG$ on a $C^*$-algebra $\clA$ is a non-degenerate $*$-homomorphism $\alpha :\clA\rightarrow M(c_{0}(\hG)\otimes \clA)$ satisfying $(\hD\otimes \mbox{id}_{\clA})\circ \alpha= (\mbox{id}_{c_{0}(\hG)}\otimes \alpha)\circ \alpha$ and such that $\alpha (\clA)(c_{0}(\hG)\otimes 1)$ is dense in $c_{0}(\hG)\otimes \clA$. Here we assume $\clA\subset B(\clH)$.

The reduced crossed product of $\clA$ by the action $\alpha $ of the discrete quantum group $\hG$ is the $C^*$-subalgebra of  $B(L^2(\cG)\otimes \clH)$ generated by $\alpha (\clA)$ and $C^{*}_{r}(\hG)\otimes 1$. The reduced crossed product algebra is denoted by $\clA\rtimes_{\alpha, r }\hG$. We omit the definition of the \textit{universal crossed product} (denoted by $\clA\rtimes_{\alpha}\hG$) as for amenable actions of $\hG$, the reduced and the universal algebras coincide.

\section{Proofs of the Main Results}

\subsection*{Notation} We follow the notations as given below (which are in contrast with those in \cite{SZ10}): The co-multiplication of the discrete quantum group $(\hG,\hD)$ is given by the multiplicative unitary $\bbW$, i.e. $\hD(\cdot)=\bbW^* (1\otimes \cdot)\bbW$. 

\bigskip

The following definition is motivated from \cite{SV14}. Recall that, we have $l^{\infty}(\hG) \subset B(L^2(\cG))$ and $L^{\infty}(\cG)=C^{*}_{r}(\hG)^{''}$. In what follows, we denote the center of a von-Neumann algebra $M$ by $Z(M)$.

\begin{definition}[Quantum weak expectation] \label{qwe}
Let $\hG$ be a discrete quantum group. A quantum weak expectation is a unital completely positive map $$\Phi :B(L^2(\cG))\rightarrow C^{*}_{r}(\hG)^{''}$$ such that $\Phi (x)= x$  for every $x\in C^{*}_{r}(\hG)$ and satisfying $\Phi (\ell^{\infty}(\hG))\subset Z(L^{\infty}(\cG))$.
\end{definition}

Now we give the proof of Theorem \ref{dqgwe} below.

\subsection*{Proof of Theorem \ref{dqgwe}}
First, assume that $C^*_r(\hG)$ admits a quantum weak expectation
$\Phi :B(L^2(\cG))\rightarrow C^*_r(\hG)''$ such that
$\Phi (a)=a$ for every $a\in C^*_r(\hG)$ satisfying $\Phi (\ell^{\infty}{(\hG)})\subset Z(L^{\infty}{(\cG)})$.
Let $h$ be the Haar state on the compact quantum group $\cG$.
Consider the restriction of the state $h\circ \Phi$ on $B(L^2(\cG))$ to $\ell^{\infty}{(\hG)}$, and let
$$m:= h\circ \Phi|_{\ell^{\infty}(\hG)}.$$
We will show that $m$ is a left invariant mean on $\ell^{\infty}{(\hG)}$.
It is sufficient to prove that 
$m((\hW \otimes \mbox{id})\hD(x))= \hW(1)m(x)$
for every $\hW\in \ell^{1}(\hG)_+$, $x\in \ell^{\infty}(\hG)_1^+$. To show this, we explicitly use the structure of $\ell^{\infty}(\hG)$, the structure of the associated unitary $\mathbb{W}$ and the assumption that the set $\irrG$ is countable.

Fix $\hW\in \ell^{1}(\hG)_+$ and $x\in \ell^{\infty}(\hG)_1^+$. Recall that, for the discrete quantum group $\hG=(\ell^{\infty}(\hG), \hD)$, there exists a family of central, mutually orthogonal, finite dimensional projections $(z_{\alpha})_{\alpha\in \irrG}$ such that $$\ell^{\infty}(\hG)=\prod_{\alpha\in \irrG}\ell^{\infty}(\hG)z_{\alpha}=\prod_{\alpha\in \irrG} B(H_\alpha)$$ and $1_{B(L^2(\cG))}=\sum_{\alpha\in \irrG}z_{\alpha}$ (in SOT). Let $\hW_\alpha(\cdot):=\hW(z_\alpha \cdot)$. Then $\hW=\sum \hW_\alpha$ and $\hW_\alpha(y)=\mbox{tr}_\alpha (\beta_\alpha z_\alpha y)$ for any $y\in \ell^{\infty}(\hG)$, some $\beta_\alpha\geq 0$ and with $\mbox{tr}_\alpha$ being the normalized trace in $B(H_\alpha)$. We have $$\|\hW\|= \sum \mbox{tr}_\alpha \beta_\alpha < \infty.$$Next, let $\epsilon >0$. Choose a finite subset $F\subset\subset \irrG$ such that the functional $$\hW_F:=\sum_{\alpha \in F} \hW_\alpha$$satisfies $\|\hW-\hW_F\|< \frac{\epsilon}{2}$. Note that, $\hW_F \leq \hW$ and $\mbox{supp } \hW_F=z_F:=\sum_{\alpha\in F} z_\alpha$. The multiplicative unitary $\mathbb{W}$ has the form $$\bigoplus_\alpha \sum_{i_\alpha, j_\alpha} e^\alpha_{i_\alpha, j_\alpha} \otimes a^\alpha_{i_\alpha, j_\alpha},$$ where $e^\alpha_{i_\alpha, j_\alpha}$ are appropriate matrix units and $a^\alpha_{i_\alpha, j_\alpha} \in C^*_r(\hG)$. Now, consider the expression $(\hW_F\otimes \mbox{id})\hD(x)$. Using the fact that $(z_{\alpha}\otimes 1)\bbW=\bbW(z_{\alpha}\otimes 1)$ and the expression of the co-multiplication $\hD$, it is clear that $$(\hW_F\otimes \mbox{id})\hD(x)=\sum_{\alpha, \alpha^\prime \in F} \hW_F( (e^\alpha_{i_\alpha, j_\alpha})^* (e^{\alpha^\prime}_{i_{\alpha^\prime}, j_{\alpha^\prime}})) (a^\alpha_{i_\alpha, j_\alpha})^* x (a^{\alpha^\prime}_{i_{\alpha^\prime}, j_{\alpha^\prime}}).$$ Observe that, the sum on the right hand side of the equation above is finite. Consider the extension of the operator $\hD$ to $B(L^2(\cG))$ by the recipe $$T\longmapsto \mathbb{W}^* (1\otimes T)\mathbb{W}$$ for $T\in B(L^2(\cG))$. By slight abuse of notation, we continue to denote this extension by $\hD$, where the appropriate definition is applicable depending on the argument of this operator. The fact that $\Phi$ is a weak expectation, the expression $(\hW_F\otimes \mbox{id})\hD(x)$ has finitely many terms in the summation formula and using the extended definition of $\hD$ we have $$\Phi\circ (\hW_F \otimes \mbox{id})\hD(x)= (\hW_F\otimes \mbox{id})\hD(\Phi(x)),$$since $C^*_r(\hG)$ lies in the multiplicative domain of $\Phi$ by virtue of the bimodule property of $\Phi$ from Definition \ref{qwe}. This equality leads to the following computation:
\begin{align*}
m((\hW_F \otimes \mbox{id})\hD(x))&=h\circ \Phi ((\hW_F \otimes \mbox{id})(\mathbb{W}^{*}(1\otimes x)\mathbb{W}))\\
&=h((\hW_F \otimes \mbox{id})(\mathbb{W}^{*}(1\otimes \Phi(x))\mathbb{W}))\\
&=h((\hW_F \otimes \mbox{id})(1\otimes \Phi(x)))\\
&=h(\hW_F (1)\Phi(x))\\
&=\hW_F (1)(h\circ \Phi(x))\\
&=\hW_F (1)m(x),
\end{align*}
where the third equality is due to the condition of quantum weak expectation i.e., $\Phi(\ell^{\infty}(\hG))\subset Z(L^{\infty}(\cG))$. To finish off this direction of the proof we have the estimate:
\begin{align*}
|m((\hW &\otimes \mbox{id})\hD(x))- \hW(1)m(x)| \leq \\
&|m(((\hW-\hW_F) \otimes \mbox{id})\hD(x))|+ |m((\hW_F \otimes \mbox{id})\hD(x)) - \hW(1)m(x)|\\
\leq &\|m\| \|((\hW-\hW_F) \otimes \mbox{id})\hD(x))\| + |\hW_F (1)m(x) - \hW(1)m(x)|\\
\leq &\|\hW-\hW_F\| \|\hD(x)\| + |(\hW_F (1)- \hW(1))m(x)|\\
\leq &\frac{\epsilon}{2} + \frac{\epsilon}{2} = \epsilon,
\end{align*}
as $x\in \ell^{\infty}(\hG)_1^+$ and the first part of the third inequality is a well-known fact. Since $\epsilon$ was arbitrary we have the desired result.

Conversly, assume that $\hG$ is amenable. Then by \cite{SV14}, $\cG$  is \textit{quantum injective} (c.f. \cite{SV14, B} for definition) i.e., there exists a \textit{conditional expectation} $E : B(L^2(\cG))\rightarrow L^{\infty}(\cG)$ such that $E(\ell^{\infty}(\hG))\subset Z(L^{\infty}(\cG))$ which by definition is also a quantum weak expectation. \qed                                                

\begin{remark}
In the proof of Theorem \ref{dqgwe} given above, observe that the non-traciality of the Haar functional of a general compact quantum group is of no consequence owing to the nature of the quantum weak expectation. Indeed, the Haar functional may as well be replaced by any positive functional in $L^{\infty}(\cG)^*$ to yield the same result. Further, if $\cG$ is of the Kac-type (i.e. the Haar functional is tracial), then the assumption of the existence of an \textit{ordinary} weak expectation is enough to establish the existence of an invariant mean on $\ell^{\infty}(\hG)$. However, by \cite{SV14}, even in this case, the existence of a quantum weak expectation is guaranteed.    
\end{remark}

\begin{remark}
The results of \cite{SV14} hold in the general context of locally compact quantum groups. As mentioned earlier, our definition of the quantum weak expectation is motivated from the notion of \textit{quantum injectivity}. However, due to the lack of an explicit description of the von-Neumann algebra structure of $L^{\infty}(\hG)$, Peter-Weyl duality and the structure of the associated unitary $\bbW$ in the locally compact quantum group setting; our proof of Theorem \ref{dqgwe} holds true only in the context described in the statement of the theorem.
\end{remark}

Next, we give the proof of Theorem \ref{cr}.

\subsection*{Proof of Theorem \ref{cr}}  Let $\hG$ be an amenable discrete quantum group and $\clA\subset B(\clH_{\clA})$ be the universal representation of $\clA$ on its universal Hilbert space $\clH_\clA$. Assume that $\clA$ has the weak expectation property. Since we have $\hG$ to be amenable, henceforth we simply denote the unique crossed product by $\clA\rtimes_{\alpha}\hG$ as the reduced and universal crossed products are the same.

The faithful action $\alpha$ yields a canonical copy of $\clA$ inside the crossed product algebra by $\alpha (\clA)\cong \alpha (\clA)(1_{C^*_r(\hG)}\otimes 1) \subset \clA\rtimes_{\alpha }\hG$ \cite [Remark 2.5.]{SZ10}. So we have the following faithfully represented inclusions:
$$\alpha(\clA) (1_{C^*_r(\hG)}\otimes 1) \subset \clA\rtimes_{\alpha }\hG\subset B(L^2(\cG) \otimes \clH_\clA).$$

Denote (for ease of notation) the universal representation of $\clA\rtimes_{\alpha }\hG$ by $\pi$ on the universal Hilbert space $\clH_u$ of $\clA\rtimes_{\alpha }\hG$. Since $\clA$ has the weak expectation property, there exists a ucp map $\Phi_\clA :B(\clH_\clA)\rightarrow \clA^{\prime\prime}$ such that $\Phi_\clA (a)=a$ for all $a\in \clA$. Next, consider the representation $\pi_0 : \clA \rightarrow B(\clH_u)$ defined by $\pi_0 (a)=\pi (\alpha (a))$ and denote its normal extension to $\clA^{\prime\prime}$ by $$\pi^{\prime\prime}: \clA^{\prime\prime} \rightarrow B(\clH_u).$$ We have the von-Neumann algebraic inclusion $\pi^{\prime\prime}(\clA^{\prime\prime}) \subset \pi({\clA\rtimes_{\alpha }\hG})^{\prime\prime}$. 

For the discrete quantum group $\hG=(c_0(\hG), \hD)$, there exists a family of central, mutually orthogonal, finite dimensional projections $(z_{i})_{i\in I}$ such that $c_0(\hG)=\bigoplus_{i\in \irrG}c_0(\hG)z_{i}$ and $1_{B(L^2(\cG))}=\sum_{i\in \irrG}z_{i}$ (in SOT). If $F$ is a finite non-empty subset of $\irrG$, then $z_F=\sum_{i\in F}z_{i}\in Z(c_0(\hG))\subset B(L^2(\cG))$ is a finite rank orthogonal projection and $z_{F} B(L^2(\cG))z_{F}=B(z_{F}L^2(\cG))$.

Let $\{e_p\}_{p=1}^{m}$ be an orthonormal basis of the finite dimensional Hilbert space $H_{F}:=z_{F}L^2(\cG)$. In what follows, we adopt notations which are similar to those used in the proof of \cite [Theorem 3.1.]{SZ10} if not the same.

Define a ucp map $\phi_{F}:B(L^2(\cG) \otimes \clH_\clA)\rightarrow B(\clH_{F})\otimes B(\clH_\clA)$ by $$\phi_{F}(y)=(z_{F}\otimes 1)y(z_{F}\otimes 1)$$  for $y\in B(L^2(\cG) \otimes \clH_\clA)$. Note that, for $y=\alpha(a)(x\otimes 1)\in\AG$ one has $$\phi_F(\alpha(a)(x\otimes 1))=(z_F\otimes 1)\alpha(a)(z_F x z_F\otimes 1).$$

Next, we consider the unitary operator in $B(L^2(\cG)\otimes H_u)$ given by: $$\widetilde{\bbW}:=(\id \otimes \pi) (\bbW \otimes 1) $$ and for $\xi \in H_F$, define the operator $\widetilde{V}_\xi\in B(H_F\otimes H_u, H_u)$ given as a row matrix by:$$\widetilde{V}_\xi= [\tilde{\lambda}(\hW_{\xi, e_1}) \cdots  \tilde{\lambda}(\hW_{\xi, e_p})]$$where $\tilde{\lambda}(\hW_{\xi, e_j})=(\hW_{\xi, e_j} \otimes \id) \widetilde{\bbW}$; $j=1,\ldots, p$. From \cite[Theorem 3.1.]{SZ10} we know that $\widetilde{V}_\xi$ is a co-isometry. Let $\id_F$ denote the map $\id_{B(H_F)}$.

Define a ucp map $\Theta_\xi : B(\clH_\clA \otimes L^2(\cG)) \rightarrow \pi({\clA\rtimes_{\alpha }\hG})^{\prime\prime}$ by the formula: $$\Theta_\xi := R_{\widetilde{V}_\xi} \circ (\id_F \otimes \pi^{\prime\prime})\circ (\id_F \otimes \Phi_\clA) \circ \phi_F$$ where $R_{\widetilde{V}_\xi} (T) = \widetilde{V}_\xi T \widetilde{V}_\xi^*$ for $T\in B(H_F\otimes H_u)$ and $\Phi_\clA$ is the weak expectation.

It is clear from the definition of $\Theta_\xi$ and the fact that $\tilde{\lambda}(\hW) \in \pi({\clA\rtimes_{\alpha }\hG})^{\prime\prime}$ for any $\hW \in \ell^1(\hG)$, that $\mbox{range }\Theta_\xi \subset \pi({\clA\rtimes_{\alpha }\hG})^{\prime\prime}$. Now, we compute the image of $\alpha(a)(x\otimes 1)\in\AG$ under the action of $\Theta_\xi$ below:
\begin{align*}
\Theta_\xi &(\alpha(a)(x\otimes 1)) = (R_{\widetilde{V}_\xi} \circ (\id_F \otimes \pi^{\prime\prime})) ((z_F\otimes 1)\alpha(a)(z_F x z_F\otimes 1))\\
& = R_{\widetilde{V}_\xi} ((\id_F \otimes \pi)((z_F\otimes 1 \otimes 1)((\id \otimes \alpha) \alpha(a))(z_F x z_F\otimes 1\otimes 1)))\\
& = \pi ((V_\xi \otimes 1)(z_F\otimes 1 \otimes 1)((\id \otimes \alpha) \alpha(a))(z_F x z_F\otimes 1\otimes 1)(V_\xi^* \otimes 1))\\
& = \pi (\alpha(a)(T_{\omega_\xi}(x)\otimes 1)),
\end{align*} where $V_{\xi}$ is the operator as defined in the proof of \cite[Theorem 3.1.]{SZ10} and $T_{\omega_\xi}$ is the convolution operator on $C^*_r(\hG)$. The last equality is a direct consequence of \cite[Equation (3.1)]{SZ10}.

Since $\hG$ is amenable, by Tomatsu \cite{Tomatsu06}, there exist a net of finitely supported vectors $\{\xi_j\}_j\subset L^2(\cG)$, i.e. $z_{F_j} \xi_j = \xi_j$ for all $j$ and $F_j \subset \subset I$, such that $T_{\omega_{\xi_j}}(x)\rightarrow x$ for all $x\in C^*_r(\hG)$. Let $\{\Theta_{\xi_j}\}_j$ be the net of ucp maps defined above corresponding  to the net of vectors $\{\xi_j\}_j$ and let $$\Theta:=\mbox{pt.wk.}-\lim_j \Theta_{\xi_j}.$$ Then $\Theta : B(L^2(\cG) \otimes \clH_\clA) \rightarrow \pi(\clA\rtimes_{\alpha }\hG)^{\prime\prime}$ and $\Theta|_{\clA\rtimes_{\alpha }\hG}=\pi$. By virtue of the injectivity of $B(L^2(\cG) \otimes \clH_\clA)$, we have a weak expectation $\widetilde{\Theta}$ for the universal representation $\pi$ of $\clA\rtimes_{\alpha }\hG$, which proves that $\clA\rtimes_{\alpha }\hG$ has the weak expectation property.

Conversely, suppose $\clA\rtimes_{\alpha }\hG$ has the weak expectation property. By injectivity of the map (action) $\alpha$ there exists a conditional expectation from $\clA\rtimes_{\alpha }\hG$ onto $\clA$ \cite[Lemma 2.2]{SZ10}. Therefore, $\clA$ has the weak expectation property. This concludes our proof. \qed

\section*{Acknowledgement} 
The authors gratefully thank the anonymous referee for a very careful reading of the first submitted draft of this paper, pointing out some mistakes/discrepancies in the notations used and providing several helpful suggestions which improved the overall presentation of the article. 

The first named author is supported by the second named author's Science and Engineering Board (DST, Govt. of India) grant no. ECR/2018/000016. Both authors gratefully acknowledge the generous funding from SERB.  

The authors hereby declare that they have no conflict of interest.


\end{document}